\documentclass[12pt]{article}

\begin{document}

\title{Modularity of $p$-adic Galois representations via $p$-adic 
approximations \\ {\small{\it{Dedicated to the memory of my mother Nalini
B. Khare \\ 30th August 1936--12th March 2002} }}}

\author{Chandrashekhar Khare}

\date{}

\maketitle

\newtheorem{theorem}{Theorem}
\newtheorem{lemma}{Lemma}
\newtheorem{prop}{Proposition}
\newtheorem{cor}{Corollary}
\newtheorem{conj}{Conjecture}
\newtheorem{guess}{Guess}
\newtheorem{remark}{Remark}
\newtheorem{example}{Example}
\newtheorem{conjecture}{Conjecture}
\newtheorem{definition}{Definition}
\newtheorem{quest}{Question}
\newtheorem{ack}{Acknowledgemets}
\newcommand{\rhobar}{\overline{\rho}}
\newcommand{\Sha}{{\rm III}}

\noindent{\bf Abstract:} In this short note we give a new approach
to proving modularity of $p$-adic Galois representations
using a method of $p$-adic approximations.

\section*{Modularity lifting theorem}

In the work of Wiles in [W], as completed by Taylor and Wiles in [TW], 
the modularity of many 2-dimensional $p$-adic representations
of the absolute Galois group $G_{\bf Q}$ of ${\bf Q}$ was
proven assuming that the mod $p$ reduction of the representation
was irreducible and modular. The proof was via proving
the isomorphism of certain deformation and Hecke rings. 
A more naive approach to proving the modularity
of a $p$-adic representation, say 
$$\rho:G_{\bf Q} \rightarrow GL_2({\bf Z}_p),$$
assuming that its reduction $\rhobar$ is modular,
that works directly with $\rho$ instead of fitting it into a 
family (i.e., interpreting it as a point 
in the spectrum of a deformation ring), 
and then proving modularity for the family 
as is done in [W] and [TW],
would be as follows: Starting from the assumption that $\rhobar$ is modular
prove successively that the mod $p^n$ reductions
of $\rho$ occur in the $p$-power torsion of the abelian variety $J_1(N)$
for a fixed $N$. In this note we give a proof of 
modularity lifting results in this more direct style.  
Here like in [K] we merely want to present a new method for proving
known results, and will illustrate the method by rederiving
the following special case of the results proven in
[W] and [TW]. This is not the optimal result that can be obtained
by this method: see the end of the note for the statement of some refinements.
 
\begin{theorem} (A. Wiles, R. Taylor)
  Let $\rho:G_{\bf Q} \rightarrow GL_2(W(k))$ be a continous representation,
  with $W(k)$ the Witt vectors of a finite field $k$ of residue characteristic
  $p>5$. 

  Assume that the mod $p$ reduction $\rhobar$ of $\rho$ has the
  following properties: \begin{itemize} \item ${\rm Ad}^0(\rhobar)$ is absolutely
  irreducible, \item  
  $\rhobar$ is modular.\end{itemize}

  Further assume that: \begin{itemize} \item $\rho$ is semistable at all
  primes, \item $\rho$ is of weight 2 at $p$ and Barsotti-Tate at $p$ if $\rhobar$ is
  finite, flat at $p$, \item
  and that the primes ramified in $\rho$ are finitely many and 
  not $\pm 1$ mod $p$. \end{itemize}

  Then $\rho$
  arises from $S_2(\Gamma_0(N))$ for some integer $N$.
\end{theorem}

\noindent{\bf Remarks:} 

1. The idea of such a proof was proposed in [K1], but at that time we 
could not put it into practise. In [K1] we had observed that for many $\rho$'s
(for examples the ones in the theorem), assuming $\rhobar$ is modular one 
can show that $\rho_n$ arises from $J_1(N(n))$ for a positive integer $N(n)$ that depends
on $n$. The new observation of the present note is that in many circumstances 
using this we can deduce (see Proposition 1) that $\rho_n$ arises from $J_1(N)$ for a {\it fixed} $N$.

2. By semistable at primes $q \neq p$ we simply mean that the
restriction to the inertia at $q$ should be unipotent, and at $p$ semistable of weight 2
we mean that $\rho$ at $p$ should either be {\it Barsotti-Tate}, i.e.,
arise from a $p$-divisible group, or be ordinary and
the restriction to the inertia at $p$ should be of the form $\left(\matrix{\varepsilon&*\cr
                                                                             0&1}\right)$
with $\varepsilon$ the $p$-adic cyclotomic character. Note that the determinant of such a $\rho$ is
$\varepsilon$.

\section*{Proof} The rest of the paper will be occupied with the proof of this theorem.
The proof will have 2 steps: We first prove that the reduction mod $p^n$ of $\rho$, 
$\rho_n:G_{\bf Q} \rightarrow GL_2(W_n(k))$ with $W_n(k)$ the Witt
vectors of length $n$ of $k$, arises in the $p$-power torsion of $J_0(Q_n)$
where $Q_n$ is the (square-free) product of a finite set of primes that depends on $n$. 
For this we use the Ramakrishna-lifts or $R$-lifts of [R], the determination
of their limit points in Theorem 1 of [K1]. In the second step, 
we deduce from the first step that
$\rho_n$ arises from $J_0(N)$ for a 
positive integer $N$ independent of $n$ (see Proposition 1 below).
From this we will deduce the theorem easily.

\subsection*{Step 1}

Let $S$ be the set of ramification of $\rhobar$. We repeat the following lemma
from [K1] and its proof for convenience. To explain the notation used,
by {\it $R$-primes} we
mean primes $q$ that are not $\pm 1$ mod $p$, unramified for the residual representation 
$\rhobar$ and for which $\rhobar({\rm Frob}_q)$ has eigenvalues with ratio $q^{\pm 1}$.
We say that a finite set of {\it $R$-primes} $Q$ is {\it auxiliary} if certain maps on $H^1$ and
$H^2$, namely $H^1(G_{S \cup Q},{\rm Ad}^0(\rhobar))
\rightarrow \oplus _{v \in S \cup Q}H^1(G_v,{\rm Ad}^0(\rhobar))/{\cal
N}_v$
and $H^2(G_{S \cup Q},{\rm Ad}^0(\rhobar))
\rightarrow \oplus _{v \in S \cup Q} H^2(G_v,{\rm Ad}^0(\rhobar))$
considered in [R] (and which we
refer to for the notation used: recall that ${\cal N}_v$ for $v \in Q$
is the mod $p$ cotangent space of a smooth quotient of the local
deformation ring at $v$ which parametrises {\it special} lifts) are
isomorphisms.
These isomorphisms result in the fact that there is a {\it unique} 
lift $\rho_{S \cup Q}^{Q-new}:G_{\bf Q} \rightarrow GL_2(W(k))$ which is semistable
of weight 2, unramified outside $S \cup Q$, minimally ramified at
primes in $S$, with determinant $\varepsilon$, and {\it special} at primes in $Q$
(by {\it special} at $q$ we mean that the restriction 
locally at the decomposition group $D_q$ at $q$ should up to twist be of the form
$\left(\matrix{\varepsilon&*\cr
                0&1}\right)$: we use this definition even 
for representations into $GL_2(R)$ with $R$ a $W(k)$-algebra).

\begin{lemma}
 Let $Q_n'$ be any finite set primes that includes the primes of ramification of $\rho_n$, 
such that $Q_n' \backslash S$ contains only $R$-primes and
such that $\rho_n|_{D_q}$ is special for $q \in Q_n' \backslash S$.
Then there exists a finite set of primes $Q_n$ that contains $Q_n'$,
such that $\rho_n|_{D_q}$ is special for $q \in Q_n \backslash S$, $Q_n \backslash S$ contains only $R$-primes
and $Q_n \backslash S$ is auxiliary. 
\end{lemma}

\noindent{\bf Proof:} We use [R] and Lemma 8 of [KR] (that latter
being a certain mutual disjointness result for field extensions cut out by
$\rho_n$ and extensions cut out by  elements of 
$H^1(G_{\bf Q},{\rm Ad}^0(\rhobar))$ and $H^1(G_{\bf Q},{\rm Ad}^0(\rhobar)^*)$) 
to construct an auxiliary set of primes $V_n$
such that $\rho_n|_{D_q}$ is special for $q \in V_n$. Then
as $Q_n' \backslash S$ contains only $R$-primes, it follows from Proposition 1.6 of [W]
that the kernel and cokernel of the map $$H^1(G_{S \cup V_n \cup Q_n'},{\rm Ad}^0(\rhobar))
\rightarrow \oplus _{v \in S \cup V_n \cup Q_n'}H^1(G_v,{\rm Ad}^0(\rhobar))/{\cal N}_v$$
have the same cardinality, as the domain and range have the same
cardinality. Then using Proposition 10 of [R],
or Lemma 1.2 of [T], and Lemma 8 of [KR], 
we can augment the set $S \cup V_n \cup Q_n'$ to get a set $Q_n$ as in the statement of the lemma.
 
\vspace{3mm}

\noindent{\bf Remark:} We can choose $Q_n'$ as in the lemma such that
$Q_n'$ is independent of $n$ (as $\rho$ is ramified at only finitely
many primes). But the set $Q_n$ that the lemma
produces depends much on $n$, and can be chosen to be independent of
$n$ only if $\rho$ is itself a $R$-lift. Further note that just like the
auxiliary primes in [TW], the sets $Q_n$ have no coherence property in
general.

\vspace{3mm}

We choose a finite set of primes $Q_n'$ as in Lemma 1 and use the lemma
to complete $Q'_n$ to a set $Q_n$ such that $Q_n \backslash S$ is auxiliary and
$\rho_n|_{D_q}$ is special for $q \in Q_n \backslash S$.
Then we claim $\rho_{Q_n}^{Q_n \backslash S-new} \equiv \rho$ mod $p^n$.
The claim is true, as the set $Q_n \backslash S$ being auxiliary,
there is a {\it unique} representation $G_{\bf Q} \rightarrow
GL_2(W(k)/(p^n))$ (with determinant $\varepsilon$)
that is unramified outside $Q_n$, minimal at $S$ and
special at primes of $Q_n \backslash S$. (It is of vital importance that
$\rho$ is $GL_2(W(k))$-valued as otherwise we would not be 
able to invoke the disjointness results
that are used in the proof of Lemma 1 (Lemma 8 of [KR]).)

\vspace{3mm}

Because of the uniqueness alluded to above, it follows 
from the level-raising results of [DT] (see Theorem 1 of [K]) that
$\rho_{Q_n}^{Q_n \backslash S-new}$
arises from $J_0(Q_n)$ (where abusively we denote by $Q_n$ the product of
primes in $Q_n$), and hence because of the congruence 
$\rho_{Q_n}^{Q_n \backslash S-new} \equiv \rho$ mod $p^n$, we deduce that
$\rho_n$ arises from (i.e., is isomorphic as a $G_{\bf Q}$-module to a
submodule of) the $p$-power torsion of
$J_0(Q_n)$ and for primes $r$ prime to $Q_n$, $T_r$ acts on $\rho_n$
via ${\rm tr}(\rho({\rm Frob}_r))$. 

\vspace{3mm}

\noindent{\bf Remark:} Its worth noting that in the proof above the
property of auxiliary sets $Q$ that gets used is that $R_{S \cup
Q}^{Q-new}$ is a (possibly finite) quotient of $W(k)$. Thus it 
is the uniqueness of lifts with given local properties rather than
their existence which is crucial for the work here (as also in [K]).

\subsection*{Step 2}

Let $W_n$ be the subset of
$Q_n$ at which $\rho_n$ is unramified (note that the set $Q_n
\backslash W_n$
is independent of $n$ for $n>>0$ as $\rho$ is {\it finitely} ramified).
Then we have the proposition:

\begin{prop}
  $\rho_n$ arises from the $W_n$-old subvariety of $J_0(Q_n)$, and furthermore
  all the Hecke operators $T_r$, for $r$ a prime not dividing $Q_n$, act on $\rho_n$ by
  ${\rm tr}(\rho_n({\rm Frob}_r))$.
\end{prop}

\noindent{\bf Proof:} This is a simple application of Mazur's principle
(see Section 8 of [Ri]). The principle relies 
on the fact that on torsion points of Jacobians 
of modular curves with semistable reduction at a prime $q$, which 
are unramified at $q$ and which reduce to lie in the 
``toric part'' of the reduction mod $q$ of
these Jacobians, the Frobenius action is constrained. Namely,
on the ``toric part'' the Frobenius ${\rm Frob}_q$ acts by $-w_qq$
where $w_q$ is the Atkin-Lehner involution.
We flesh this out this below.

Consider a
prime $q \in W_n$. Then decompose $\rho_n|_{D_q}$ (which is unramified by hypothesis) 
into $W(k)/p^n \oplus W(k)/p^n$ where
on the first copy ${\rm Frob}_q$ acts by a scalar that is not $\pm q$:
this is possible as $q^2$ is not 1 mod $p$ and $\rho_n$ is special at
$q$. Let $e_n$ be a generator for the first summand. 
We would like to prove 
that $\rho_n$ occurs in the $q$-old subvariety of
$J_0(Q_n)$. Note that using irreducibility of $\rhobar$, 
Burnside's lemma gives that $\rhobar(k[G_{\bf Q}])=M_2(k)$ and hence
by Nakayama's lemma $\rho_n(W_n(k)[G_{\bf Q}])=M_2(W_n(k))$.
Thus using the fact that the $q$-old subvariety is
stable under the Galois and Hecke action, the fact that $\rho_n$ 
occurs in the $q$-old subvariety of
$J_0(Q_n)$ is implied by the claim that
$e_n$ is contained in the $q$-old subvariety of $J_0(Q_n)$. Let $\cal J$ be the N\'eron model 
at $q$ of $J_0(Q_n)$. Note that
as $\rho_n$ is unramified at $q$ it maps injectively to ${\cal
J}_{/{\bf F}_q}(\overline{{\bf F}_q})$ under
the reduction map. Now if the claim were false, as the group of connected components
$\cal J$ is Eisenstein (see loc. cit.), we would deduce that
the reduction of $e_n$ in ${\cal J}^0_{/{\bf F}_q}(\overline{\bf F}_q)$
maps non-trivially (and hence its image 
has order divisible by $p$) to the $\overline{{\bf F}_q}$-points of
the torus which is the quotient of ${\cal J}_{/{\bf F}_q}^0$ 
by the image of the $q$-old subvariety (in characteristic $q$).
But as we recalled above, it is well known (see loc. cit.) that ${\rm Frob}_q$ 
acts on the $\overline{\bf F}_q$-valued points of this toric quotient
(isogenous 
to the torus $T$ of ${\cal J}_{{\bf
F}_q}^0$, the latter being a semiabelian variety that is 
an extension of $J_0({Q_n \over q})_{/{\bf F}_q}^2$ by $T$) 
by $-w_qq$ which gives the contradiction that $q^2$ is 1 mod $p$. 
This contradiction proves the claim. Now taking another prime $q' \in W_n$
and working within the $q$-old subvariety of $J_0(Q_n)$, by the same argument we
see that $\rho_n$ occurs in the $\{q,q'\}$-old subvariety of $J_0(Q_n)$, and
eventually that $\rho_n$ occurs in the $W_n$-old subvariety of $J_0(Q_n)$.
Furthermore by inspection the last part of the proposition is also clear.



\begin{cor}
  $\rho_n$ arises from $J_0(N)$ for a fixed integer $N$ (independent of $n$).
\end{cor}

\noindent{\bf Proof:} This follows from the proposition from the
irreducibility of $\rhobar$ and the fact that the kernel of
the degeneracy map from several copies of $J_0({{Q_n} \over {W_n}})$
to $J_0(Q_n)$ is Eisenstein (see [Ri1]), where by abuse we
denote by $W_n$ the product of the primes in $W_n$.

\subsection*{$\rho$ is modular}

There are various ways by which we can deduce 
the theorem from the corollary and we choose the following.
From the corollary, we see that $\rho_n$ occurs in $J_0(N)$ for a fixed $N$, such that
for all $r$ not a divisor of any of the $Q_n$'s, the Hecke operator $T_r$ acts on $\rho_n$ via the
reduction mod $p^n$ of the scalar (that is really an endomorphism!) 
${\rm tr}(\rho({\rm Frob}_r))$. 
By Proposition 1 of [K1], which proves that the density of 
primes that ramify in any element of a converging sequence of residually absolutely irreducible $p$-adic
Galois representations is 0, we see that the set consisting of all prime
divisors of all of the $Q_n$'s has density 0. From this we deduce that for a density 1 set of primes $r$ there is a 
normalised newform $f$ in $S_2(\Gamma_0(N))$ such that $a_r(f)={\rm tr}(\rho({\rm Frob}_r))$: this is because
the intersection of the kernels of $T_r-{\rm tr}(\rho({\rm Frob}_r))$ 
(for such primes $r$) 
acting on the $p$-primary torsion group of $J_0(N)$ has unbounded
exponent (this shows that we could work systematically with the
$W_n$-old subvariety of $J_0(Q_n)$ and thus only use the proposition above 
instead of the corollary). 
Hence by the Cebotarev density theorem
$\rho$ arises from $f$ concluding the proof of the theorem.

\vspace{3mm}

\noindent{\bf Remarks:} 

1. In [K], the $R$-lifts of [R] 
were used to give new proofs of modularity theorems
that did not use TW systems but nevertheless generally relied on the set-up in [W]
of comparing deformation and Hecke rings and the numerical criterion for
isomorphisms of complete intersections of [W].

2. By cutting the Jacobians we work with into pieces according to the action of 
the Atkin-Lehner involutions we can make the proof of Proposition 1 work when
the prime $q$ is not 1 mod $p$.

3. While the theorem makes many hypotheses on $\rho$ the 
one that is most vital for the method to work, besides the usual hypotheses
that are present in [W] et al, is that $\rho$ is $W(k)$-valued. The others
are of a more technical nature and perhaps some of them could be eased
after some further
work. (Note that for a $GL_2$-type abelian variety the corresponding $p$-adic
representation is defined over unramified extensions of ${\bf Q}_p$
for almost all primes $\wp$.) As far as the ramification hypotheses 
on $\rho$ go these can certainly be relaxed: for instance using the
refinement of the result of [R] in [T],
and the methods here, we can obtain the following result.

\begin{theorem}
  Let $\rhobar:G_{\bf Q} \rightarrow GL_2(k)$ be a continous, odd representation, with $k$ a finite
  field of characteristic bigger than 3, such that $Ad^0(\rhobar)$
  is irreducible. Assume that $\rhobar$ is modular, and at $p$
  is up to twist neither the trivial representation nor unramified
  with image of
  order divisible by $p$. 

  Let $\rho:G_{\bf Q}
  \rightarrow GL_2(W(k))$ 
  be a continuous lift of $\rhobar$ that has the following properties:
  \begin{itemize} \item
  $\rho$ is minimally ramified 
  at the primes of ramification ${\rm Ram}(\rhobar)$ of $\rhobar$,
  \item $\rho$ is of weight 2 at $p$,
  and  Barsotti-Tate at $p$ if $\rhobar$ is
  finite, flat at $p$, \item the set of primes ${\rm Ram}(\rho)$ ramified in $\rho$ is
  finite, \item $\rho$ is semistable at all
  the primes of ${\rm Ram}(\rho) \backslash {\rm Ram}(\rhobar)$, \item
  and for the primes $q$ in ${\rm Ram}(\rho)$ that are 
  not in ${\rm Ram}(\rhobar)$, 
  $\rhobar|_{D_q}$ is not a scalar. \end{itemize} 
  
  Then $\rho$ arises from a newform of weight 2.
\end{theorem}

It will be of interest to have a less restrictive theorem accessible
by the methods of this paper, for instance be able to treat (many)
3-adic representations. Conditions ensuring minimality of
ramification of $\rho$ at ${\rm Ram}(\rhobar) \cup {p}$ seem
essential. 

\vspace{3mm}

\noindent{\bf Acknowledgements:} I would like to thank Srinath Baba
for a conversation in May, 2002 in Montreal which prompted me to think again
of the idea of [K1] after a lapse of 2 years.

\section*{References}

\noindent [DT] Diamond, F., Taylor, R., 
{\it Lifting modular mod $l$ representations},
Duke Math. J. 74 (1994), no. 2, 253--269.

\vspace{3mm}

\noindent [K] Khare, C., {\it On isomorphisms between deformation rings and
Hecke rings}, preprint available at {\sf http://www.math.utah.edu/\~{ }shekhar/papers.html}

\vspace{3mm}

\noindent [K1] Khare, C., {\it Limits of residually irreducible
$p$-adic Galois representations},to appear in Proc. of the AMS,
available at {\sf http://www.math.utah.edu/\~{ }shekhar/papers.html}

\vspace{3mm}

\noindent [KR] Khare, C., Ramakrishna, R., {\it
Finiteness of Selmer groups and deformation rings}, preprint available at
{\sf http://www.math.utah.edu/\~{ }shekhar/papers.html}

\vspace{3mm}

\noindent [Ri] Ribet, K., {\it Report on mod $\ell$
representations of ${\rm Gal}(\overline{\bf Q}/{\bf Q})$}, in Motives, Proc.
Sympos. Pure Math. 55, Part 2 (1994), 639--676.

\vspace{3mm}

\noindent [Ri1] Ribet, K., {\it Congruence relations between modular
forms}, Proc. International Cong. of Math. (1983), pp 503--514.

\vspace{3mm}

\noindent [R] Ramakrishna, R., {\it Deforming Galois representations and the
conjectures of Serre and Fontaine-Mazur}, to appear in Annals of Math.

\vspace{3mm}

\noindent [T] Taylor, R., {\it On icosahedral Artin
representations II}, to appear in American J of Math.

\vspace{3mm}

\noindent[TW] Taylor, R., Wiles, A., 
{\it Ring-theoretic properties of certain Hecke
algebras}, Ann. of Math. (2) 141 (1995), 553--572. 

\vspace{3mm}

\noindent [W] Wiles, A., {\it Modular elliptic curves and 
Fermat's last theorem}, Ann. of Math. 141 (1995), 443--551. 

\vspace{3mm}

\noindent {\it Address of the author}: Dept of Math, Univ of Utah, 155 S 1400 E,
Salt lake City, UT 84112. e-mail address: shekhar@math.utah.edu
 
\noindent School of Mathematics, 
TIFR, Homi Bhabha Road, Mumbai 400 005, INDIA. e-mail addresses: shekhar@math.tifr.res.in

\end{document}